\documentclass[12pt]{article}
\usepackage{wasysym,amssymb,indentfirst,cite,bibspacing,mathrsfs,amsthm}
\begin{document}
\newtheorem{theorem}{THEOREM}[section]
\newtheorem{lemma}[theorem]{LEMMA}
\newtheorem{proposition}[theorem]{PROPOSITION}
\newtheorem{example}[theorem]{EXAMPLE}
\newtheorem{definition}[theorem]{DEFINITION}
\newtheorem{corollary}[theorem]{COROLLARY}
\baselineskip 17pt
\title{\textbf{FINITE GROUPS IN WHICH SS-PERMUTABILITY IS A TRANSITIVE RELATION}\thanks{Research is supported by a NNSF grant of China (grant \#11371335) and Research Fund for the Doctoral Program of Higher Education of China (Grant 20113402110036).}}
\author{X. Y. CHEN and W. B. GUO\thanks{Corresponding author.}\\
{\small Wu Wen-Tsun Key Laboratory of Mathematics, and} \\ {\small  Department of Mathematics, University of Science and Technology of China,}\\ {\small Hefei, 230026, P. R. China}\\
 {\small e-mails: jelly@mail.ustc.edu.cn, $\,$wbguo@ustc.edu.cn}
}
\date{}
\maketitle

\textbf{Abstract.} A subgroup $H$ of a finite group $G$ is said to be SS-permutable in $G$ if $H$ has a supplement $K$ in $G$ such that $H$ permutes with every Sylow subgroup of $K$. A finite group $G$ is called an SST-group if SS-permutability is a transitive relation on the set of all subgroups of $G$. The structure of SST-groups is investigated in this paper.

\renewcommand{\thefootnote}{\empty}

\section{Introduction}
Throughout this paper, all groups considered are finite. For a group $G$, $\pi(G)$ denotes the set of prime divisors of $|G|$, $|G|_p$ denotes the order of a Sylow $p$-subgroup of $G$, and $Syl_p(G)$ denotes the set of all Sylow $p$-subgroups of $G$.\par
A subgroup $H$ of a group $G$ is said to be \textit{permutable} (resp. \textit{S-permutable}) in $G$ if $H$ permutes with all the subgroups (resp. Sylow subgroups) of $G$. A group $G$ is called a \textit{T-group} (resp. \textit{PT-group}, \textit{PST-group}) if normality (resp. permutability, S-permutability) is a transitive relation, that is, if $H$ and $K$ are subgroups of $G$ such that $H$ is normal (resp. permutable, S-permutable) in $K$ and $K$ is normal (resp. permutable, S-permutable) in $G$, then $H$ is normal (resp. permutable, S-permutable) in $G$.\par
By Kegel's result \cite{Keg}, a group $G$ is a PST-group if and only if every subnormal subgroup of $G$ is S-permutable in $G$. Recall that the nilpotent residual of a group $G$ is the intersection of all normal subgroups $N$ of $G$ such that $G/N$ is nilpotent. Agrawal \cite{Agr} showed that a group $G$ is a solvable PST-group if and only if the nilpotent residual $L$ of $G$ is a normal abelian Hall subgroup of $G$ upon which $G$ acts by conjugation as power automorphisms. In particular, a solvable PST-group is supersolvable. The structure of PST-groups has been investigated by many authors, see for example \cite{Ale,Bal1,Bal2,Bei,Bal3,Bal}.\par
A subgroup $H$ of a group $G$ is said to be \textit{semipermutable} (resp. \textit{S-semipermutable}) \cite{Che} in $G$ if $H$ permutes with every subgroup (resp. Sylow subgroup) $X$ of $G$ such that $(|H|,|X|)=1$. A group $G$ is called a \textit{BT-group} (resp. an \textit{SBT-group}) \cite{Wan} if semipermutability (resp. S-semipermutability) is a transitive relation. The following result is useful in the sequel, which was established by Wang et al. in \cite{Wan}.\par
\begin{theorem}[\cite{Wan}, Theorem 3.1]\label{Theorem 1.1}Let $G$ be a group. Then the following statements are equivalent:\par
$(1)$ $G$ is a solvable BT-group.\par
$(2)$ $G$ is a solvable SBT-group.\par
$(3)$ Every subgroup of $G$ is semipermutable in $G$.\par
$(4)$ Every subgroup of $G$ is S-semipermutable in $G$.\par
$(5)$ Every subgroup of $G$ of prime power order is semipermutable in $G$.\par
$(6)$ Every subgroup of $G$ of prime power order is S-semipermutable in $G$.\par
$(7)$ $G$ is a solvable PST-group, and if $p$ and $q$ are distinct primes not dividing the order of $L$, where $L$ is the nilpotent residual of $G$, then $[G_p,G_q]=1$ for every $G_p\in Syl_p(G)$ and $G_q\in Syl_q(G)$.
\end{theorem}
Recall that a subgroup $H$ of a group $G$ is said to be \textit{SS-permutable} (or \textit{SS-quasinormal}) in $G$ if $H$ has a supplement $K$ in $G$ such that $H$ permutes with every Sylow subgroup of $K$. In this case, $K$ is called an SS-permutable supplement of $H$ in $G$. This important embedding property was introduced by Li et al. in \cite{Li}, and they called these subgroups SS-quasinormal subgroups. In this paper, we use the term SS-permutable subgroups for such subgroups. It is easy to see that every SS-permutable subgroup of a group $G$ is S-semipermutable in $G$ (see below Lemma \ref{Lemma 2.1}(4)). However, the converse does not hold in general as the following example shows.\par
\begin{example}\label{Example 1.2}\textup{Let $G=\langle x,y,z,w \,|\, x^3=y^3=z^2=w^2=1, [x,y]=[z,w]=1, x^z=x^{-1}, y^w=y^{-1}, x^w=x^{-1}, y^z=xy\rangle$ and $H=\langle y,w \rangle$. Obviously, $H$ is S-semipermutable in $G$. Now assume that $H$ is SS-permutable in $G$ with an SS-permutable supplement $K$. Let $K_2$ be a Sylow 2-subgroup of $K$. Then by definition, $H$ permutes with $K_2$. It is clear that $\langle y \rangle \unlhd H\unlhd HK_2$, and so $HK_2\leq N_G(\langle y \rangle)$. This implies that $\langle y \rangle\unlhd G$, which is contrary to our assumption. Therefore, $H$ is not SS-permutable in $G$.}\end{example}
We say that a subgroup $H$ of a group $G$ is \textit{NSS-permutable} in $G$ if $H$ has a normal supplement $K$ in $G$ such that $H$ permutes with every Sylow subgroup of $K$. In this case, $K$ is called an NSS-permutable supplement of $H$ in $G$. Note that every NSS-permutable subgroup of a group $G$ is SS-permutable in $G$, but the next example illustrates that the converse is not true.\par
\begin{example}\label{Example 1.3}\textup{Let $G=A_5$ and $H=A_4$, where $A_5$ and $A_4$ denote the alternating group of degree 5 and 4, respectively. Clearly, $K\in Syl_5(G)$ is an SS-permutable supplement of $H$ in $G$. This shows that $H$ is SS-permutable in $G$. If $H$ is NSS-permutable in $G$, then $G$ is the only NSS-permutable supplement of $H$ in $G$. This induces that $H$ is S-permutable in $G$, and so $H$ is subnormal in $G$ by \cite[Theorem 1]{Keg}, which is impossible. Hence $H$ is not NSS-permutable in $G$.}\end{example}
We say that a group $G$ is an \textit{SST-group} (resp. \textit{NSST-group}) if SS-permutability (resp. NSS-permutability) is a transitive relation. Motivated by the above-mentioned results on PST-groups, BT-groups and SBT-groups, we naturally consider to investigate the structure of SST-groups and NSST-groups. Thus, the purpose in this paper is to establish several characterizations of SST-groups and NSST-groups.\par
Note that a group is said to be an \textit{SC-group} if all its chief factors are simple. Clearly, the class of all supersolvable groups coincides with the class of all solvable SC-groups. We begin with:\par
\begin{theorem}\label{Theorem A}Let $G$ be an SST-group $($resp. NSST-group$)$. Then $G$ is an SC-group.\end{theorem}
\begin{theorem}\label{Theorem B}Let $G$ be a group. Then the following statements are equivalent:\par
$(1)$ $G$ is solvable and every subnormal subgroup of $G$ is SS-permutable in $G$.\par
$(2)$ $G$ is solvable and every subnormal subgroup of $G$ is NSS-permutable in $G$.\par
$(3)$ Every subgroup of $F^*(G)$ is SS-permutable in $G$.\par
$(4)$ Every subgroup of $F^*(G)$ is NSS-permutable in $G$.\par
$(5)$ $G$ is a solvable PST-group.\end{theorem}
In Theorem \ref{Theorem B} above, $F^*(G)$ denotes the generalized Fitting subgroup of $G$, that is, the product of all normal quasinilpotent subgroups of $G$. The readers may refer to \cite[Chapter X]{Hup} for details.\par
\begin{theorem}\label{Theorem C}Let $G$ be a group. Then the following statements are equivalent:\par
$(1)$ Whenever $H\leq K$ are two $p$-subgroups of $G$ with $p\in\pi(G)$, $H$ is SS-permutable in $N_G(K)$.\par
$(2)$ Whenever $H\leq K$ are two $p$-subgroups of $G$ with $p\in\pi(G)$, $H$ is NSS-permutable in $N_G(K)$.\par
$(3)$ $G$ is a solvable PST-group.\end{theorem}
Our main result is the following:\par
\begin{theorem}\label{Theorem D}Let $G$ be a solvable group. Then the following statements are equivalent:\par
$(1)$ $G$ is an SST-group.\par
$(2)$ $G$ is an NSST-group.\par
$(3)$ Every subgroup of $G$ is SS-permutable in $G$.\par
$(4)$ Every subgroup of $G$ is NSS-permutable in $G$.\par
$(5)$ Every subgroup of $G$ of prime power order is SS-permutable in $G$.\par
$(6)$ Every subgroup of $G$ of prime power order is NSS-permutable in $G$.\par
$(7)$ Every cyclic subgroup of $G$ of prime power order is SS-permutable in $G$.\par
$(8)$ Every cyclic subgroup of $G$ of prime power order is NSS-permutable in $G$.\end{theorem}
\begin{corollary}\label{Corollary 1.4}A solvable SST-group $G$ is a BT-group.\end{corollary}
The following example illustrates that a solvable BT-group is not necessarily a solvable SST-group.\par
\begin{example}\label{Example 1.5}\textup{Let $G=\langle x,y \,|\,x^{5}=y^4=1,\,x^y=x^{2} \rangle$. Note that the nilpotent residual of $G$ is $\langle x\rangle\in Syl_{5}(G)$. By Theorem \ref{Theorem 1.1}, $G$ is a solvable BT-group. Put $H=\langle y \rangle$ and $L=\langle y^2 \rangle$. Suppose that $L$ is SS-permutable in $G$. Then $G$ is the unique supplement of $L$ in $G$. It follows that $L$ is S-permutable in $G$, and thus $L\leq O_2(G)$ by \cite[Theorem 1]{Keg}. This implies that either $O_2(G)=H$ or $O_2(G)=L$. It is easy to verify that $y^{x}=yx^{-1}$ and ${(y^2)}^{x}=y^2x^2$, and so $H\ntrianglelefteq G$ and $L\ntrianglelefteq G$, a contradiction. Hence $L$ is not SS-permutable in $G$. Then by Theorem \ref{Theorem D}, $G$ is not a solvable SST-group.}\end{example}
\begin{corollary}\label{Corollary 1.6}The class of all solvable SST-groups is closed under taking subgroups and epimorphic images.\end{corollary}
Recall that a subgroup $H$ of a group $G$ is called \textit{abnormal} in $G$ if $g\in \langle H,H^g\rangle$ for all $g\in G$.\par
\begin{corollary}\label{Corollary 1.7}Let $G$ be a solvable group. Then the following statements are equivalent:\par
$(1)$ $G$ is an SST-group.\par
$(2)$ Every subgroup of $G$ is either SS-permutable or abnormal in $G$.\par
$(3)$ Every subgroup of $G$ is either NSS-permutable or abnormal in $G$.\end{corollary}
Let $G$ be a solvable group and let $\pi(G)=\{p_1,\ldots, p_k\}$. Suppose that $Q_i$ is a Hall ${p_i}'$-subgroup of $G$ for $1\leq i\leq k$. Then the set $\{Q_1,\ldots,Q_k\}$ is called a \textit{Sylow system} of $G$, and the subgroup $H=\bigcap_{i=1}^{k} N_G(Q_i)$ is called a \textit{system normalizer} of $G$.\par
Suppose that $G$ is a solvable PST-group. Let $L$ be the nilpotent residual of $G$ and $D$ be a system normalizer of $G$. Then $L$ is a normal abelian Hall subgroup of $G$ upon which $D$ acts by conjugation as power automorphisms. Hence by \cite[(9.2.7)]{Rob}, $D$ is a Hall subgroup of $G$ and $G=L\rtimes D$. Note that all the system normalizers of $G$ are conjugate in $G$.\par
Recall that a group $G$ is said to be \textit{complemented} if every subgroup of $G$ has a complement in $G$. In \cite{Hal}, Hall proved that the class of all complemented groups is exactly the class of all supersolvable groups with elementary abelian Sylow subgroups.\par
Next we give another two characterizations of solvable SST-groups.\par
\begin{theorem}\label{Theorem E}Let $G$ be a solvable SST-group with the nilpotent residual $L$ and a system normalizer $D$. Then the Frattini group $\Phi(G)=\Phi(L)\times \Phi(D)$. Furthermore, $G/\Phi(G)$ is a complemented group.\end{theorem}
\begin{theorem}\label{Theorem F}Let $G$ be a solvable BT-group with the nilpotent residual $L$. Then the following statements are equivalent:\par
$(1)$ $G$ is an SST-group.\par
$(2)$ For every $p$-subgroup $P$ of $G$ with $p\in \pi(G)\backslash \pi(L)$, $G$ has a $p$-subgroup $K_p$ such that $PK_p\in Syl_p(G)$ and $[P,\langle {K_p}^L\rangle]\leq O_p(G)$.\par
$(3)$ For every $p$-element $x$ of $G$ with $p\in \pi(G)\backslash \pi(L)$, $G$ has a $p$-subgroup $K_p$ such that $\langle x \rangle K_p\in Syl_p(G)$ and $[x,\langle {K_p}^L\rangle]\leq O_p(G)$.\end{theorem}
It was established in \cite{Chi} that $G_1\times G_2 \times \cdots \times G_n$ ($n\geq 2$) is a solvable PST-group if and only if $G_i$ is a solvable PST-group with the nilpotent residual $L_i$ for $1\leq i\leq n$ and $(|L_i|,|G_j|)=1$ for $1\leq i,j \leq n$ and $i\neq j$. Now we continue the study for solvable BT-groups and solvable SST-groups.\par
\begin{theorem}\label{Theorem G}$G_1\times G_2 \times \cdots \times G_n$ $($$n\geq 2$$)$ is a solvable BT-group if and only if $G_i$ is a solvable BT-group with the nilpotent residual $L_i$ for $1\leq i\leq n$ and $(|L_i|,|G_j|)=1$ for $1\leq i,j \leq n$ and $i\neq j$.\end{theorem}
\begin{theorem}\label{Theorem H}$G_1\times G_2 \times \cdots \times G_n$ $($$n\geq 2$$)$ is a solvable SST-group if $G_i$ is a solvable SST-group and $(|G_i|,|G_j|)=1$ for $1\leq i,j \leq n$ and $i\neq j$.\end{theorem}
The next example shows that Theorem \ref{Theorem G} does not hold for solvable SST-groups.\par
\begin{example}\label{Example 1.8}\textup{Let $G=\langle x,y,z,w \,|\, x^3=y^5=z^2=w^2=1, [x,y]=[x,w]=[z,y]=[z,w]=1, x^z=x^{-1}, y^w=y^{-1}\rangle$, $G_1=\langle x,z\rangle$ and $G_2=\langle y,w\rangle$. It is clear that $G_1$ and $G_2$ are solvable SST-groups and $L_1=\langle x\rangle$ and $L_2=\langle y\rangle$ are the nilpotent residuals of $G_1$ and $G_2$, respectively. Put $H=\langle zw\rangle$. It is easy to verify that $G$, $\langle x,y,z\rangle$ and $\langle x,y,w\rangle$ are not SS-permutable supplements of $H$ in $G$. Hence $H$ is not SS-permutable in $G$. Then by Theorem \ref{Theorem D}, $G$ is not a solvable SST-group.}\end{example}
In \cite{Bei1}, J. C. Beidleman et al. proved that a group $G$ is a solvable PST-group if and only if $G$ has a normal subgroup $N$ such that $N$ and $G/N''$ are solvable PST-groups. Analogously, we have the following result.\par
\begin{theorem}\label{Theorem I}A group $G$ is a solvable SST-group $($resp. BT-group$)$ if and only if $G$ has a normal subgroup $N$ such that $N$ is a solvable PST-group and $G/N''$ is a solvable SST-group $($resp. BT-group$)$.\end{theorem}
\section{Preliminaries}
\begin{lemma}\label{Lemma 2.1}Suppose that a subgroup $H$ of a group $G$ is SS-permutable $($resp. NSS-permutable$)$ in $G$ with an SS-permutable supplement $($resp. NSS-permutable supplement$)$ $K$, $L\leq G$ and $N\unlhd G$. Then:\par
$(1)$ If $H\leq L$, then $H$ is SS-permutable $($resp. NSS-permutable$)$ in $L$.\par
$(2)$ $HN/N$ is SS-permutable $($resp. NSS-permutable$)$ in $G/N$.\par
$(3)$ If $N\leq L$ and $L/N$ is SS-permutable $($resp. NSS-permutable$)$ in $G/N$, then $L$ is SS-permutable $($resp. NSS-permutable$)$ in $G$.\par
$(4)$ $H$ is S-semipermutable in $G$.\par
$(5)$ If $H\leq F(G)$, then $H$ is S-permutable in $G$.\par
$(6)$ Every conjugate of $K$ in $G$ is an SS-permutable supplement $($resp. NSS-permutable supplement$)$ of $H$ in $G$.\par
$(7)$ If $N$ is nilpotent, then $NK$ is an SS-permutable supplement $($resp. NSS-permutable supplement$)$ of $H$ in $G$.\par
$(8)$ If $H$ is a $p$-subgroup of $G$ with $p\in \pi(G)$, then $HK_p\in Syl_p(G)$ for every $K_p\in Syl_p(K)$.\end{lemma}
\begin{proof}[\textup{PROOF}] We only prove the statements for SS-permutable subgroups. For NSS-permutable subgroups, the statements can be handled similarly.\par
(1)-(3) See \cite[Lemma 2.1]{Li}.\par
(4) By definition, $H$ permutes with every Sylow subgroup of $K$. Let $X$ be a Sylow subgroup of $G$ such that $(|H|,|X|)=1$. Then there exists an element $h\in H$ such that $X^h\leq K$. It follows that $HX^h=X^hH$, and so $HX=XH$. Hence $H$ is S-semipermutable in $G$.\par
(5) See \cite[Lemma 2.2]{Li}.\par
(6) The statement is obvious.\par
(7) For every prime $p\in\pi(K)$, there exist $N_p\in Syl_p(N)$ and $K_p\in Syl_p(K)$ such that $N_pK_p\in Syl_p(NK)$. Note that $N_p\unlhd G$. Then for any Sylow $p$-subgroup $P$ of $NK$, $NK$ has an element $g=kh$ for some $k\in K$ and $h\in H$ such that $P={(N_pK_p)}^g={(N_p{K_p}^k)}^h$. Since $K$ is an SS-permutable supplement of $H$ in $G$, we have that $H(N_p{K_p}^k)=(N_p{K_p}^k)H$, and so $HP=H{(N_p{K_p}^k)}^h={(N_p{K_p}^k)}^hH=PH$. Therefore, $NK$ is an SS-permutable supplement of $H$ in $G$.\par
(8) Evidently, $K$ has an element $k$ such that $H{K_p}^k\in Syl_p(G)$. Note that $H\cap K$ is S-permutable in $K$, and so it is subnormal in $K$ by \cite[Theorem 1]{Keg}. Hence $H\cap K\leq O_p(K)$. This implies that $H\cap K=H\cap K_p=H\cap {K_p}^k$. It follows that $|G|_p=|H{K_p}^k|=|H||{K_p}^k|/|H\cap {K_p}^k|=|H||K_p|/|H\cap K_p|=|HK_p|$. Thus $HK_p\in Syl_p(G)$.\end{proof}
\begin{lemma}[\cite{Bal4}, Lemma 2.1.3]\label{Lemma 2.2}Let $p$ be a prime and $N$ a normal $p$-subgroup of a group $G$. Then all subgroups of $N$ are S-permutable in $G$ if and only if all chief factors of $G$ below $N$ are cyclic and $G$-isomorphic when regarded as modules over $G$.\end{lemma}
Recall that a subgroup $H$ of a group $G$ is said to be \textit{$\tau$-quasinormal} \cite{Luk1} in $G$ if $HG_p=G_pH$ for every $G_p\in Syl_p(G)$ such that $(|H|,p)=1$ and $(|H|,|{G_p}^G|)\neq 1$. It is easy to see that every SS-permutable subgroup of $G$ is $\tau$-quasinormal in $G$ by Lemma \ref{Lemma 2.1}(4).\par
\begin{lemma}[\cite{Luk}, Theorem 1.2]\label{Lemma 2.3}Let $G$ be a group. Then every subgroup of $F^*(G)$ is $\tau$-quasinormal in $G$ if and only if $G$ is a solvable PST-group.\end{lemma}
\begin{lemma}\label{Lemma 2.4}A solvable group $G$ is a BT-group if and only if every cyclic subgroup of $G$ of prime power order is S-semipermutable in $G$.\end{lemma}
\begin{proof}[\textup{PROOF}] By Theorem \ref{Theorem 1.1}, we only need to prove the sufficiency. Assume that every cyclic subgroup of $G$ of prime power order is S-semipermutable in $G$. Let $H$ be a $p$-subgroup of $G$ with $p\in \pi(G)$. Then for every element $h\in H$ and every Sylow $q$-subgroup $Q$ of $G$, we have that $\langle h \rangle Q=Q\langle h \rangle$. Consequently, we have that $HQ=QH$, and so $H$ is S-semipermutable in $G$. This shows that every subgroup of $G$ of prime power order is S-semipermutable in $G$. Hence by Theorem \ref{Theorem 1.1}, $G$ is a BT-group.\end{proof}
\begin{lemma}\label{Lemma 2.5}Let $G$ be a nilpotent-by-abelian group and $H\leq G$. Then $H$ is SS-permutable in $G$ if and only if $H$ is NSS-permutable in $G$.\end{lemma}
\begin{proof}[\textup{PROOF}] The sufficiency is clear. We only need to prove the necessity. Since $G$ is nilpotent-by-abelian, $G/F(G)$ is abelian. Suppose that $H$ is SS-permutable in $G$ with an SS-permutable supplement $K$. Then by Lemma \ref{Lemma 2.1}(7), $F(G)K$ is also an SS-permutable supplement of $H$ in $G$. As $F(G)K\unlhd G$, we have that $F(G)K$ is an NSS-permutable supplement of $H$ in $G$. Hence $H$ is NSS-permutable in $G$.\end{proof}
\begin{lemma}\label{Lemma 2.6}Let $T$ and $S$ are SS-permutable $($resp. NSS-permutable$)$ in a solvable group $G$ such that $(|T|,|S|)=1$. Then $\langle T,S\rangle$ is SS-permutable $($resp. NSS-permutable$)$ in $G$. The solvability assumption is not necessary if $T$ and $S$ are subgroups of prime power order.\end{lemma}
\begin{proof}[\textup{PROOF}] We only prove the lemma for SS-permutable subgroups. Let $K_1$ and $K_2$ be SS-permutable supplements of $T$ and $S$ in $G$, respectively. Note that $G$ is solvable. By Lemma \ref{Lemma 2.1}(6), without loss of generality, we may assume that $S\leq K_1$ and $T\leq K_2$. Then $TS(K_1\cap K_2)=TK_1=G$. This means that $K_1\cap K_2$ is a supplement of $\langle T,S\rangle$ in $G$. For any Sylow $p$-subgroup $P$ of $K_1\cap K_2$ with $p\in \pi(K_1\cap K_2)$, there exist a Sylow $p$-subgroup $K_{1p}$ of $K_1$ and a Sylow $p$-subgroup $K_{2p}$ of $K_2$ such that $P=K_{1p}\cap K_2=K_1\cap K_{2p}$. Since $T(K_{1p}\cap K_2)=(K_{1p}\cap K_2)T$ and $S(K_1\cap K_{2p})=(K_1\cap K_{2p})S$, we have that $\langle T,S\rangle P=P\langle T,S\rangle$. This shows that $K_1\cap K_2$ is an SS-permutable supplement of $\langle T,S\rangle$ in $G$. Thus $\langle T,S\rangle$ is SS-permutable in $G$.\end{proof}
\begin{lemma}\label{Lemma 2.7}Let $G$ be a solvable PST-group with the nilpotent residual $L$ and a system normalizer $D$. Then:\par
$(1)$ $D'\unlhd G$.\par
$(2)$ If $H$ is an SS-permutable $p$-subgroup of $G$ with an SS-permutable supplement $K$, where $p\in \pi(G)$, then $[H,K_p]\leq O_p(G)$ for every $K_p\in Syl_p(K)$.\end{lemma}
\begin{proof}[\textup{PROOF}] (1) Note that for every element $l\in L$, $\langle l \rangle$ is normal in $G$. Then $G/C_G(l)$ is abelian, and so $D'\leq G'\leq C_G(l)$. It follows that $D'\leq C_G(L)$. Hence $D'\unlhd G$ for $G=L\rtimes D$.\par
(2) If $p\in \pi(L)$, then it is evident. Now assume that $p\in \pi(G)\backslash \pi(L)$. Then for every $K_p\in Syl_p(K)$, $G$ has a Sylow $p$-subgroup $G_p$ such that $G_p=HK_p$ by Lemma \ref{Lemma 2.1}(8). It follows from (1) that ${G_p}'\unlhd G$. Hence $[H,K_p]\leq {G_p}'\leq O_p(G)$.\end{proof}
\begin{lemma}\label{Lemma 2.8}Let $G$ be a solvable BT-group with the nilpotent residual $L$. Suppose that $H$ is a $p$-subgroup of $G$ with $p\in \pi(G)\backslash \pi(L)$. Then $H$ is SS-permutable in $G$ if $G$ has a $p$-subgroup $K_p$ such that $HK_p\in Syl_p(G)$ and $[H,\langle {K_p}^L\rangle]\leq O_p(G)$.\end{lemma}
\begin{proof}[\textup{PROOF}] Let $D$ be a system normalizer of $G$ containing $HK_p$ and $D_{p'}$ be the Hall $p'$-subgroup of $D$. Then we claim that $K=O_p(G)K_pD_{p'}L$ is an SS-permutable supplement of $H$ in $G$. In fact, it is clear that $G=HK$. By Theorem \ref{Theorem 1.1}, if $q\in \pi(K)$ and $q\neq p$, then $H$ permutes with every Sylow $q$-subgroup of $K$. Note that $D_{p'}\leq C_G(O_p(G)K_p)$. Then for every element $g\in K$, we have that $[H,{(O_p(G)K_p)}^g]\leq O_p(G)$. This deduces that $H{(O_p(G)K_p)}^g={(O_p(G)K_p)}^gH$. Hence $H$ permutes with every Sylow $p$-subgroup of $K$. Thus the claim holds, and consequently, $H$ is SS-permutable in $G$.\end{proof}

\section{Proof of the Theorems}
\begin{proof}[{\textbf{\textup{PROOF OF THEOREM \textup{\ref{Theorem A}}.}}}] We only prove the theorem for SST-groups. For NSST-groups, the theorem can be proved similarly. Let $N$ be a minimal normal subgroup of $G$. By Lemma \ref{Lemma 2.1}(2) and (3), $G/N$ is also an SST-group. By induction, $G/N$ is an SC-group. Now we only need to prove that $N$ is simple. Since $G$ is an SST-group, all subnormal subgroups of $G$ are SS-permutable in $G$. It follows that all subnormal subgroups of $N$ are SS-permutable in $G$. If $N$ is abelian, then $N$ is a normal $p$-subgroup of $G$ for some $p\in \pi(G)$. By Lemma \ref{Lemma 2.1}(5), all subgroups of $N$ are S-permutable in $G$. Hence by Lemma \ref{Lemma 2.2}, $N$ is cyclic, and so $N$ is simple. Now assume that $N$ is non-abelian. Then $N=N_1\times N_2\times \cdots \times N_n$, where $N_i$ is isomorphic to a non-abelian simple group $S$ for $1\leq i\leq n$. Note that $N_1$ is SS-permutable in $G$. Let $K$ be an SS-permutable supplement of $N_1$ in $G$. Then $N_1K_p=K_pN_1$ for every $K_p\in Syl_p(K)$. Since ${N_1}^g\unlhd N$ for every element $g\in G$, we have that $\langle{N_1}^{K_p}\rangle \unlhd N$. By \cite[Chapter A, Proposition 4.13(b)]{Doe}, $\langle{N_1}^{K_p}\rangle$ is a direct product of a subset of $N_i$ ($1\leq i\leq n$). Obviously, $\langle{N_1}^{K_p}\rangle /N_1\leq N_1K_p/N_1$ is a $p$-group. This deduces that $\langle{N_1}^{K_p}\rangle=N_1$, and thereby $K_p\leq N_G(N_1)$. It follows that $K\leq N_G(N_1)$, and thus $N_1\unlhd G$. Hence $N=N_1$ is simple. This completes the proof of Theorem \ref{Theorem A}.\end{proof}
\begin{proof}[{\textbf{\textup{PROOF OF THEOREM \textup{\ref{Theorem B}}.}}}] Note that every NSS-permutable subgroup of a group $G$ is SS-permutable in $G$. Hence (2) implies (1), and (4) implies (3). Assume that $G$ is a solvable PST-group. Then every subnormal subgroup of $G$ is S-permutable, and so NSS-permutable in $G$. Therefore, (5) implies (2).\par
Now we show that (3) implies (5). Suppose that every subgroup of $F^*(G)$ is SS-permutable in $G$. Then every subgroup of $F^*(G)$ is $\tau$-quasinormal in $G$. By Lemma \ref{Lemma 2.3}, $G$ is a solvable PST-group, and thus (5) holds. Finally, we prove that (1) implies (4). Assume that $G$ is solvable and every subnormal subgroup of $G$ is SS-permutable in $G$. Then $F^*(G)=F(G)$ by \cite[Chapter X, Corollary 13.7(d)]{Hup}. Therefore, every subgroup of $F(G)$ is SS-permutable in $G$. By Lemma \ref{Lemma 2.1}(5), every subgroup of $F(G)$ is S-permutable, and so NSS-permutable in $G$. Hence (4) follows. This ends the proof of Theorem \ref{Theorem B}.\end{proof}
\begin{proof}[{\textbf{\textup{PROOF OF THEOREM \textup{\ref{Theorem C}}.}}}] (2) implies (1) is clear. Assume that (1) holds. By Lemma \ref{Lemma 2.1}(5), whenever $H\leq K$ are two $p$-subgroups of $G$ with $p\in\pi(G)$, $H$ is S-permutable in $N_G(K)$. It follows from \cite[Theorem 4]{Bal1} that $G$ is a solvable PST-group, and so (3) follows. By \cite[Theorem 4]{Bal1} again, we also see that (3) implies (2).\end{proof}
\begin{proof}[{\textbf{\textup{PROOF OF THEOREM \textup{\ref{Theorem D}}.}}}] Firstly we show that (1) is equivalent to (3). Clearly, (3) implies (1). Now consider that $G$ is an SST-group. Then every subnormal subgroup of $G$ is SS-permutable in $G$. By Theorem \ref{Theorem B}, $G$ is a PST-group. Let $L$ be the nilpotent residual of $G$. Note that every subgroup of $L$ is normal in $G$. Then for any subgroup $H$ of $G$, $H$ is SS-permutable in $HL$. Since $HL$ is subnormal in $G$, we have that $HL$ is SS-permutable in $G$. This implies $H$ is SS-permutable in $G$. Hence (1) implies (3), and thus (1) and (3) are equivalent. With a similar argument as above, we get that (2) is equivalent to (4).\par
Next we shall prove that (5) is equivalent to (7). Clearly, (5) implies (7). Now suppose that every cyclic subgroup of $G$ of prime power order is SS-permutable in $G$. Then $G$ is a BT-group by Lemma \ref{Lemma 2.1}(4) and Lemma \ref{Lemma 2.4}. Let $L$ be the nilpotent residual of $G$ and $H$ be a $p$-subgroup of $G$ with $p\in \pi(G)$. If either $H\leq L$ or $H$ is cyclic, then there is nothing to prove. Hence we may suppose that $p\in \pi(G)\backslash \pi(L)$ and $H$ is not cyclic. If for every element $h\in H$, $\langle h \rangle$ is S-permutable in $G$, then $H$ is S-permutable in $G$ by \cite[Corollary 1]{Sch}, and so $H$ is SS-permutable in $G$. We may, therefore, assume that there exists an element $x\in H$ such that $\langle x \rangle$ is not S-permutable in $G$. By hypothesis, $\langle x \rangle$ is SS-permutable in $G$. Let $K$ be an SS-permutable supplement of $\langle x \rangle$ in $G$. By Lemma \ref{Lemma 2.1}(7), we may let $O_p(G)\leq K$. If $H\cap K=H$, then $\langle x \rangle \leq H\leq K$. This implies that $K=G$ and $\langle x \rangle$ is S-permutable in $G$, a contradiction. Hence $H\cap K<H$. By induction, $H\cap K$ is SS-permutable in $G$. Let $T$ be an SS-permutable supplement of $H\cap K$ in $G$. By Lemma \ref{Lemma 2.1}(7), we may also let $O_p(G)\leq T$. Now we claim that $K\cap T$ is an SS-permutable supplement of $H$ in $G$. Clearly, $H(K\cap T)=\langle x \rangle(H\cap K)(K\cap T)=G$. For every Sylow $q$-subgroup $Q$ of $K\cap T$ with $q\in \pi(K\cap T)$ and $p\neq q$, since $G$ is a BT-group, we have that $H$ permutes with $Q$ by Theorem \ref{Theorem 1.1}. Let $P$ be a Sylow $p$-subgroup of $K\cap T$. Then there exist a Sylow $p$-subgroup $K_p$ of $K$ and a Sylow $p$-subgroup $T_p$ of $T$ such that $P=K_p\cap T=K\cap T_p$. By Lemma \ref{Lemma 2.7}(2), $[\langle x \rangle,K_p]\leq O_p(G)$ and $[H\cap K,T_p]\leq O_p(G)$. Since $O_p(G)\leq K_p\cap T$, we have that $\langle x \rangle P=\langle x \rangle (K_p\cap T)=(K_p\cap T)\langle x \rangle=P\langle x \rangle$. Similarly, we can obtain that $(H\cap K)P=P(H\cap K)$. Therefore, $HP=\langle x \rangle(H\cap K)P=P\langle x \rangle(H\cap K)=PH$, and thus the claim follows. Then $H$ is SS-permutable in $G$, and so (5) holds. With a similar discussion as above, we get that (6) is equivalent to (8).\par
Now we show that (3) implies (4). Assume that every subgroup of $G$ is SS-permutable in $G$. By Theorem \ref{Theorem B}, $G$ is a PST-group, and so $G$ is supersolvable. It follows from Lemma \ref{Lemma 2.5} that every subgroup of $G$ is NSS-permutable in $G$. Hence (4) holds. Evidently, (4) implies (6), and (6) implies (5). Finally, assume that every subgroup of $G$ of prime power order is SS-permutable in $G$. By Lemma \ref{Lemma 2.6}, it is easy to see that every subgroup of $G$ is SS-permutable in $G$. This shows that (5) implies (3). The proof is thus completed.\end{proof}
\begin{proof}[{\textbf{\textup{PROOF OF COROLLARY \textup{\ref{Corollary 1.4}}.}}}] The corollary directly follows from Theorem \ref{Theorem D}, Lemma \ref{Lemma 2.1}(4) and Theorem \ref{Theorem 1.1}.\end{proof}
\begin{proof}[{\textbf{\textup{PROOF OF COROLLARY \textup{\ref{Corollary 1.6}}.}}}] Suppose that $G$ is a solvable SST-group. By Theorem \ref{Theorem D} and Lemma \ref{Lemma 2.1}(1), every subgroup of $G$ is a solvable SST-group. Moreover, by Lemma \ref{Lemma 2.1}(2) and (3), every epimorphic image of $G$ is also a solvable SST-group.\end{proof}
\begin{proof}[{\textbf{\textup{PROOF OF COROLLARY \textup{\ref{Corollary 1.7}}.}}}] By Theorem \ref{Theorem D}, we have that (1) implies (3). Clearly, (3) implies (2). Then we only need to prove that (2) implies (1). Suppose that every subgroup of $G$ is either SS-permutable or abnormal in $G$. Then by Lemma \ref{Lemma 2.1}(4) and \cite[Lemma 1]{Zha}, $G$ is supersolvable. Let $H$ be a $p$-subgroup of $G$ with $p\in \pi(G)$. If $p$ is not the smallest prime divisor of $|G|$, then it is easy to see that $H$ is not abnormal in $G$. This induces that $H$ is SS-permutable in $G$. We may, therefore, assume that $p$ is the smallest prime divisor of $|G|$. If $H$ is not a Sylow $p$-subgroup of $G$, then $H$ is not abnormal in $G$, and thereby $H$ is SS-permutable in $G$. If $H\in Syl_p(G)$, then for every $G_q\in Syl_q(G)$ with $q\in \pi(G)$ and $p\neq q$, since $G_q$ is SS-permutable in $G$, $H$ permutes with $G_q$. Let $K$ be the Hall $p'$-subgroup of $G$. Then $K$ is an SS-permutable supplement of $H$ in $G$, and so $H$ is SS-permutable in $G$. Hence every subgroup of $G$ of prime power order is SS-permutable in $G$. By Theorem \ref{Theorem D}, $G$ is an SST-group. This shows that (2) implies (1).\end{proof}
\begin{proof}[{\textbf{\textup{PROOF OF THEOREM \textup{\ref{Theorem E}}.}}}] Note that for every maximal subgroup $L_1$ of $L$ and every maximal subgroup $D_1$ of $D$, $L_1\rtimes D$ and $L\rtimes D_1$ are maximal subgroups of $G$. If $\Phi(G)\cap L\nleq \Phi(L)$, then $L$ has a maximal subgroup $L_1$ such that $\Phi(G)\cap L\nleq L_1$. It follows that $\Phi(G)\cap L\leq L_1D\cap L\leq L_1$, a contradiction. Hence $\Phi(G)\cap L\leq \Phi(L)$. With a similar discussion as above, we get that $\Phi(G)\cap D\leq \Phi(D)$.\par
Since $L\unlhd G$, we have that $\Phi(L)\leq \Phi(G)$. Consequently, $\Phi(G)\cap L=\Phi(L)$. Now we claim that $\Phi(P)\leq \Phi(G)$ for every Sylow $p$-subgroup $P$ of $G$ with $p\in\pi(D)$. As $G$ is a solvable SST-group, $\Phi(P)$ is SS-permutable in $G$ by Theorem \ref{Theorem D}. Let $K$ be an SS-permutable supplement of $\Phi(P)$ in $G$. Then $P=P\cap \Phi(P)K=\Phi(P)(P\cap K)$. This implies that $P\leq K$, and so $K=G$. Hence $\Phi(P)$ is S-permutable in $G$. By \cite[Lemma A]{Sch}, $O^p(G)\leq N_G(\Phi(P))$. It follows that $\Phi(P)\unlhd G$. By \cite[Chapter A, Theorem 9.2(d)]{Doe}, $\Phi(P)\leq \Phi(G)$, and thus the claim holds. Since $D$ is a nilpotent Hall subgroup of $G$, we have that $\Phi(D)\leq \Phi(G)$, and so $\Phi(G)\cap D=\Phi(D)$. Consequently, $\Phi(G)=(\Phi(G)\cap L)\times (\Phi(G)\cap D)=\Phi(L)\times \Phi(D)$. Note that $G$ is supersolvable and for every Sylow subgroup $P$ of $G$, $\Phi(P)\leq \Phi(G)$. This implies that every Sylow subgroup of $G/\Phi(G)$ is elementary abelian. Then by \cite[Theorem 2]{Hal}, $G/\Phi(G)$ is a complemented group.\end{proof}
\begin{proof}[{\textbf{\textup{PROOF OF THEOREM \textup{\ref{Theorem F}}.}}}] Firstly we show that (1) implies (2). Suppose that $G$ is an SST-group. Then by Theorem \ref{Theorem D}, every subgroup of $G$ of prime power order is SS-permutable in $G$. For every $p$-subgroup $P$ of $G$ with $p\in\pi(G)\backslash\pi(L)$, let $K$ be an SS-permutable supplement of $P$ in $G$ and $K_p$ be a Sylow $p$-subgroup of $K$. By Lemma \ref{Lemma 2.1}(8), we have that $PK_p\in Syl_p(G)$. It follows from Lemma \ref{Lemma 2.1}(6) and Lemma \ref{Lemma 2.7}(2) that $[H,{K_p}^l]\leq O_p(G)$ for every $l\in L$, and so $[H,\langle {K_p}^L\rangle]\leq O_p(G)$. Thus (2) follows.\par
It is obvious that (2) implies (3). Now we prove that (3) implies (1). Assume that (3) holds. Then by \cite[Chapter A, Lemma 7.3(c)]{Doe}, $[x^n,\langle {K_p}^L\rangle]={[x,\langle {K_p}^L\rangle]}^{x^{n-1}}\cdots \\{[x,\langle {K_p}^L\rangle]}^x[x,\langle {K_p}^L\rangle]\leq O_p(G)$ for all $n\in \mathbb{N}$. Hence $[\langle x \rangle,\langle {K_p}^L\rangle]\leq O_p(G)$. By Lemma \ref{Lemma 2.8}, $\langle x \rangle$ is SS-permutable in $G$. Note that every subgroup of $L$ is normal in $G$. Therefore, every cyclic subgroup of $G$ of prime power order is SS-permutable in $G$. Then by Theorem \ref{Theorem D}, $G$ is an SST-group, and thereby (1) holds. This completes the proof of Theorem \ref{Theorem F}.\end{proof}
\begin{proof}[{\textbf{\textup{PROOF OF THEOREM \textup{\ref{Theorem G}}.}}}] By Theorem \ref{Theorem 1.1} and \cite[Corollary 2.3]{Chi}, the necessity is obvious. Now we prove the sufficiency. By induction, we may assume that $n=2$. By \cite[Corollary 2.3]{Chi} again, $G_1\times G_2$ is a solvable PST-group. Evidently, the nilpotent residual of $G_1\times G_2$ is $L_1\times L_2$. Let $\{p,q\}\subseteq \pi(G_1\times G_2)\backslash \pi(L_1\times L_2)$ and $p\neq q$. For every Sylow $p$-subgroup $G_p$ of $G_1\times G_2$, there exist a Sylow $p$-subgroup $G_{1p}$ of $G_1$ and a Sylow $p$-subgroup $G_{2p}$ of $G_2$ such that $G_p=G_{1p}\times G_{2p}$. Similarly, for every Sylow $q$-subgroup $G_q$ of $G_1\times G_2$, there exist a Sylow $q$-subgroup $G_{1q}$ of $G_1$ and a Sylow $q$-subgroup $G_{2q}$ of $G_2$ such that $G_q=G_{1q}\times G_{2q}$. Since $G_1$ and $G_2$ are solvable BT-groups, by Theorem \ref{Theorem 1.1}, we have that $[G_{1p},G_{1q}]=1$ and $[G_{2p},G_{2q}]=1$. It follows that $[G_p,G_q]=[G_{1p}G_{2p},G_{1q}G_{2q}]=1$. Therefore, by Theorem \ref{Theorem 1.1} again, $G_1\times G_2$ is a solvable BT-group.\end{proof}
\begin{proof}[{\textbf{\textup{PROOF OF THEOREM \textup{\ref{Theorem H}}.}}}] By induction, we may assume that $n=2$. Let $L_1$ and $L_2$ be the nilpotent residuals of $G_1$ and $G_2$, respectively. By Corollary \ref{Corollary 1.4} and Theorem \ref{Theorem G}, $G_1\times G_2$ is a solvable BT-group. Obviously, the nilpotent residual of $G_1\times G_2$ is $L_1\times L_2$. Let $P$ be a $p$-subgroup of $G$ with $p\in \pi(G_1\times G_2)\backslash \pi(L_1\times L_2)$. Without loss of generality, we may suppose that $p\in \pi(G_1)\backslash \pi(L_1)$ for $(|G_1|,|G_2|)=1$. By Theorem \ref{Theorem F}, $G_1$ has a $p$-subgroup $K_p$ such that $PK_p\in Syl_p(G_1)$ and $[P,\langle {K_p}^{L_1}\rangle]\leq O_p(G_1)$. It follows that $G$ has a $p$-subgroup $K_p$ such that $PK_p\in Syl_p(G)$ and $[P,\langle {K_p}^{L_1\times L_2}\rangle]\leq O_p(G)$. By Theorem \ref{Theorem F} again, $G_1\times G_2$ is a solvable SST-group.\end{proof}
The following lemma is the main step in the proof of Theorem \ref{Theorem I}.\par
\begin{lemma}\label{Lemma 3.1}Let $G$ be a solvable PST-group, and let $Z_\infty(G)$ denote the hypercenter of $G$. If $G/Z_\infty(G)$ is a solvable SST-group $($resp. BT-group$)$, then $G$ is a solvable SST-group $($resp. BT-group$)$.\end{lemma}
\begin{proof}[\textup{PROOF}] Let $L$ be the nilpotent residual of $G$ and $D$ be a system normalizer of $G$. Firstly we prove the lemma for BT-groups. Suppose that $G/Z_\infty(G)$ is a solvable BT-group. By \cite[Chapter I, Theorem 5.9(b)]{Doe}, $Z_\infty(G)\leq D$. Note that $LZ_\infty(G)/Z_\infty(G)$ is the nilpotent residual of $G/Z_\infty(G)$. Let $G_p\in Syl_p(G)$ and $G_q\in Syl_q(G)$ with $\{p,q\}\subseteq \pi(D)$ and $p\neq q$. By Theorem \ref{Theorem 1.1}, $[G_p,G_q]\leq Z_\infty(G)$. We may, therefore, assume that $G_pG_qZ_\infty(G)\leq D$. This implies that $[G_p,G_q]=1$. By Theorem \ref{Theorem 1.1}, $G$ is a solvable BT-group.\par
Now we prove the lemma for SST-groups. Suppose that $G/Z_\infty(G)$ is a solvable SST-group. Then by above and Corollary \ref{Corollary 1.4}, $G$ is a solvable BT-group. Let $P$ be a $p$-subgroup of $G$ with $p\in \pi(G)\backslash \pi(L)$. By Theorem \ref{Theorem D}, $PZ_\infty(G)/Z_\infty(G)$ is SS-permutable in $G/Z_\infty(G)$. Let $K/Z_\infty(G)$ be an SS-permutable supplement of $PZ_\infty(G)/Z_\infty(G)$ in $G/Z_\infty(G)$. Then $K$ is a supplement of $P$ in $G$ and for every $K_p\in Syl_p(K)$, we have that $PK_pZ_\infty(G)=K_pPZ_\infty(G)$. We may, therefore, assume that $PK_pZ_\infty(G)\leq D$. Let $D_p$ be the Sylow $p$-subgroup of $D$. By Lemma \ref{Lemma 2.7}(1), $[P,K_p]\leq {D_p}'\unlhd G$. Therefore, $[P,K_p]\leq O_p(G)$. By \cite[Chapter I, Theorem 5.9(b)]{Doe}, $O_p(G)\leq \mbox{Core}_G(D)=Z_\infty(G)\leq K$, and thereby $O_p(G)\leq K_p$. This implies that $PK_p=K_pP$. By Theorem \ref{Theorem 1.1}, $P$ permutes with every Sylow $q$-subgroup of $K$ with $q\in \pi(K)$ and $q\neq p$. Hence $P$ is SS-permutable in $G$. Note that every subgroup of $L$ is normal in $G$. Therefore, every subgroup of $G$ of prime power order is SS-permutable in $G$. By Theorem \ref{Theorem D}, $G$ is a solvable SST-group.\end{proof}
\begin{proof}[{\textbf{\textup{PROOF OF THEOREM \textup{\ref{Theorem I}}.}}}] The necessity is obvious. We only need to prove the sufficiency. Let $L$ be the nilpotent residual of $G$ and $D$ be a system normalizer of $G$. By \cite[Theorem B]{Bei1}, $G$ is a solvable PST-group. Note that the class of solvable SST-groups (resp. BT-groups) is closed under taking epimorphic images. Then it is easy to see that $G/Z_\infty(G)$ satisfies the hypothesis of the theorem. If $Z_\infty(G)\neq 1$, then by induction, $G/Z_\infty(G)$ is a solvable SST-group (resp. BT-group). It follows from Lemma \ref{Lemma 3.1} that $G$ is a solvable SST-group (resp. BT-group). We may, therefore, assume that $Z_\infty(G)=1$. By Lemma \ref{Lemma 2.7}(1) and \cite[Chapter I, Theorem 5.9(b)]{Doe}, we have that $D'\leq \mbox{Core}_G(D)=Z_\infty(G)=1$. This implies that $G'=LD'=L$. Hence $N''\leq G''=1$, and so $G$ is a solvable SST-group (resp. BT-group).\end{proof}
\textbf{Acknowledgment.} The authors are grateful to the referee for his/her valuable suggestions and comments.\par

\end{document}